\documentstyle[amssymb,amsfonts]{amsart}

\def\be#1{\begin{equation} \label{#1}}
\def\bs{\begin{split}}
\def\bi{\begin{itemize}}
\def\es{\end{split}}
\def\ba{\begin{align}}
\def\bas{\begin{align*}}
\def\ea{\end{align}}
\def\eas{\end{align*}}

\def\R{{\mbox{\bf R}}}
\def\T{{\mbox{\bf T}}}

\def\dist{{\mbox{\rm dist}}}
\def\diam{{\mbox{\rm diam}}}

\def\B{{\cal B}}

\def\Z{{\mbox{\bf Z}}}
\def\q{{\mbox{\bf q}}}
\def\eps{\varepsilon}

\def\emph#1{{\it #1}}
\def\textbf#1{{\bf #1}}
\newenvironment{proof}{\noindent {\bf Proof} }{\endprf\par}
\def \endprf{\hfill  {\vrule height6pt width6pt depth0pt}\medskip}
\def\neq{{\not =}}

\theoremstyle{plain}
  \newtheorem{theorem}[subsection]{Theorem}
  \newtheorem{proposition}[subsection]{Proposition}
  \newtheorem{lemma}[subsection]{Lemma}
  \newtheorem{corollary}[subsection]{Corollary}
  
  \newtheorem{assumption}[subsection]{Assumption}

\theoremstyle{remark}

\theoremstyle{definition}
  \newtheorem{definition}[subsection]{Definition}

\include{psfig}

\begin{document}

\title[Sharp bilinear restriction for paraboloids]{ A sharp bilinear restriction estimate for paraboloids}

\author{Terence Tao}
\address{Department of Mathematics, UCLA, Los Angeles, CA 90024}
\email{tao@@math.ucla.edu}

\begin{abstract}
Recently Wolff \cite{wolff:cone} obtained a sharp $L^2$ bilinear
restriction theorem for bounded subsets of the cone in general dimension.  Here we adapt the argument of Wolff to also handle subsets of ``elliptic surfaces'' such as paraboloids.  Except for an endpoint, this answers a conjecture of Machedon and Klainerman, and also improves upon the known restriction theory for the paraboloid and sphere. 
\end{abstract}

\maketitle

\section{Introduction}

Let $n \geq 1$ be a fixed integer, and let $S$ be a smooth compact hypersurface with boundary in the space $\R \times \R^n := \{(\tau,\xi): \tau \in \R, \xi \in \R^n\}$, which we shall interpret as the spacetime frequency space.  If $0 < p, q \leq \infty$, we say that the \emph{linear adjoint restriction estimate} $\R^*_S(p \to q)$ holds if one has an estimate of the form
\be{linear}
\| \widehat{f d\sigma}\|_{L^q(\R \times \R^n)} \leq C_{p,q,S}
\|f\|_{L^p(S, d\sigma)}
\end{equation}
for all test functions $f$ on $S$, where
$$ \widehat{F}(t,x) := \int F(\tau,\xi) e^{2\pi i (t\tau + x \cdot \xi)}$$
is the spacetime Fourier transform.  The \emph{restriction problem}\footnote{Historically, the restriction problem asks for which exponents $q'$ is it true that the Fourier transform of an $L^{q'}(\R \times \R^n)$ function can be meaningfully restricted to $S$.  This is essentially the adjoint of the above problem; see \cite{stein:large} for further discussion.} for $S$ is to determine for which $p$, $q$ the estimate $\R^*_S(p \to q)$ holds.  This problem was posed by Stein \cite{stein:problem}, and is related to other outstanding problems in harmonic analysis such as the Bochner-Riesz conjecture, local smoothing conjecture, and Kakeya conjecture; see e.g. \cite{Bo}, \cite{wolff:survey}, \cite{tao:notices} for further discussion.  In one spatial dimension $n=1$, the problem is mostly solved, but in two and higher spatial dimensions the problem remains far from settled (except in special cases such as $p=2$), despite much recent progress.

It has been known for several decades that one can attack this conjecture in the special case $q=4$ by squaring both sides of the linear estimate \eqref{linear} and studying the resulting bilinear $L^2$ estimate; see e.g. \cite{feff:thesis}, \cite{sj}, etc.  Variants of this idea have also been very useful for nonlinear dispersive equations, see e.g. \cite{borg:hsd}, \cite{klainerman:nulllocal}, etc.  More recently, the same idea has been applied to more general values of $q$, see \cite{borg:cone}, \cite{tvv:bilinear}, \cite{tv:cone1}, \cite{tv:cone2}, \cite{wolff:cone}, \cite{wolff:smsub}.  More precisely, for any two smooth compact hypersurfaces  $S_1$, $S_2$ with boundary in $\R \times \R^n$, with Lebesgue measure $d\sigma_1$ and $d\sigma_2$ respectively, we say that the \emph{bilinear adjoint restriction
estimate}
$R^*_{S_1,S_2}(2 \times 2 \to q)$ holds if one has
$$ \| \widehat{f_1 d\sigma_1} \widehat{f_2 d\sigma_2} \|_{L^q(\R \times \R^n)} \leq C_{q,S_1,S_2}
\|f_1\|_{L^2(S_1, d\sigma_1)} \|f_2\|_{L^2(S_2, d\sigma_1)}.$$
for all test functions $f_1$, $f_2$ supported on $S_1$, $S_2$ respectively\footnote{One can of course place $f_1$ and $f_2$ in Lebesgue spaces other than $L^2$ (see e.g. \cite{tvv:bilinear}), but we shall not need to do so within this paper.}. 

The linear and bilinear estimates are closely related; for instance when $S_1 = S_2=S$, then $R^*_{S_1,S_2}(2 \times 2 \to q)$ is clearly equivalent to $R^*_S(2 \to 2q)$.  However, it was observed in \cite{borg:cone}, \cite{tvv:bilinear} that further estimates are available if $S_1$ and $S_2$ are not equal, and in particular if they satisfy some sort of transversality condition.  For instance, if the normals of $S_1$ and of $S_2$ are separated by at least some fixed angle $c > 0$, then one can easily obtain the bilinear estimate $R^*_{S_1, S_2}(2 \times 2 \to 2)$ by Plancherel's inequality and Cauchy-Schwarz, even in cases where the linear estimates $R^*_{S_1}(2 \to 4)$, $R^*_{S_2}(2 \to 4)$ fail.  Furthermore, these bilinear restriction estimates can then be used (via some rescaling and interpolation arguments) to obtain new linear restriction estimates; see \cite{tvv:bilinear}, \cite{tv:cone1}, \cite{wolff:cone} for some examples of this.

Two important examples of surfaces $S_1$, $S_2$ are: (a) compact, transverse subsets of the light cone 
$$\{ (\tau, \xi) \in \R \times \R^n: |\tau| = |\xi|\};$$
and (b) compact, transverse (i.e. disjoint) subsets of the paraboloid\footnote{The choice of normalization factor $-\frac{1}{2}$ may appear odd, but this is to ensure that waves of frequency $\xi$ travel at group velocity $\xi$.} 
\be{paraboloid}
S := \{ (\tau,\xi) \in \R \times \R^n: \tau = -\frac{1}{2} |\xi|^2\}.
\end{equation}
Apart from being model examples for the bilinear restriction problem, they also have direct application to nonlinear wave and Schr\"odinger equations respectively.  In 1997, Machedon and Klainerman observed that in these two cases, the estimate $R^*_{S_1,S_2}(2 \times 2 \to q)$ can only hold when $q \geq \frac{n+3}{n+1}$, and conjectured that this necessary condition was in fact sufficient (see \cite{tvv:bilinear}, \cite{wolff:cone} for further discussion).
The paraboloid \eqref{paraboloid} also serves as models for other surfaces with non-vanishing Gaussian curvature, such as the sphere; see the remarks section for further discussion.

Of the two cases (a) and (b), the cone problem was generally thought to the more difficult of the two (compare for instance \cite{borg:cone} with  \cite{borg:stein}).  It was thus a surprise when Wolff \cite{wolff:cone} established the Machedon-Klainerman conjecture for the cone in all non-endpoint cases $q > \frac{n+3}{n+1}$ (with the endpoint case being attained shortly afterward in \cite{tao:cone}). A key geometrical observation was that if one took the union of all the lines through a fixed origin $x_0$ which were normal to $S_2$, then any line normal to $S_1$ could only intersect this union in at most one point; this is ultimately due to the single vanishing principal curvature on the cone, which forces all of the above lines to be light rays.  The analogous statement for the paraboloid however is false, so one cannot directly apply Wolff's argument to case (b).  Even in two spatial dimensions $n=2$, the Machedon-Klainerman conjecture had only been verified in this case for $q > 2 - \frac{2}{17}$ (see \cite{tv:cone1}), instead of the conjectured $q \geq 2 - \frac{1}{3}$.

In this paper we adapt Wolff's argument in \cite{wolff:cone} to overcome this geometric obstruction:

\begin{theorem}\label{main}  Let $S_1$ and $S_2$ be any disjoint compact subsets of the paraboloid $S$ defined in \eqref{paraboloid}.  Then we have $R^*_{S_1, S_2}(2 \times 2 \to q)$ for any $q > \frac{n+3}{n+1}$.  In particular, the Machedon-Klainerman conjecture is true up to endpoints for the paraboloid.
\end{theorem}

By the general theory of linear and bilinear restriction theorems (see \cite{tvv:bilinear}, \cite{tv:cone1}), Theorem \ref{main} implies some new progress on the \emph{restriction conjecture} for paraboloids \cite{stein:problem}.  This conjecture asserts that $R^*(p \to q)$ holds\footnote{The numerology is shifted by one from that in \cite{stein:problem} because we are working in $\R \times \R^n$ instead of $\R^n$.} whenever $q = \frac{n+2}{n} p'$ and $q > 2(n+1)/n$, where $1/p + 1/p' = 1$; these conditions are known to be necessary.

\begin{corollary}\label{rest-conj} The restriction conjecture for paraboloids is true for $q > 2(n+3)/(n+1)$.
\end{corollary}

\begin{proof}
This follows directly from Theorem \ref{main} in this paper and Theorem 2.2 in \cite{tvv:bilinear}, together with the observation that one can freely raise the exponents $p,q$ in the estimate $R^*_{S_1, S_2}(p \times p \to q)$.  See \cite{tvv:bilinear}, \cite{tv:cone1} for more examples of this type of argument.
\end{proof}

In \cite{strichartz:restrictionquadratic} this conjecture was verified for $q \geq 2(n+2)/n$; in the special case $n=2$, the best known previous result was $q > 4 - \frac{8}{31}$ (see \cite{tv:cone1}); the above Corollary improves this to $q > 4 - \frac{2}{3}$.  (It is conjectured that this bound holds in fact for all $q > 3$).

A similar result holds for all other positively curved surfaces, such as the sphere; we discuss this in Section \ref{remarks-sec}.

Functions of the form $u := \widehat{f d\sigma}$, where $d\sigma$ is surface measure on $S$, can easily be seen to solve the \emph{free Schr\"odinger equation}
\be{schrodinger}
4\pi i u_t - \Delta u = 0.
\end{equation}
The factor $4\pi$ is an artifact of our conventions and should be ignored.
We shall call solutions to \eqref{schrodinger} \emph{free Schr\"odinger waves}.  For any free Schr\"odinger wave, the quantity $\| u(t)\|_{L^2_x(\R^n)}^2$ is an invariant of time, and shall be referred to as the \emph{total probability} $P(u)$ of the wave\footnote{This quantity plays the role of the energy for solutions to the wave equation, see \cite{tao:cone}.}.  Observe that
\be{prob}
P(\widehat{f_j d\sigma_j}) \sim \| f_j \|_2^2
\end{equation}
for any compact subset $S_j$ of $S$.

\begin{corollary}\label{bil-cor}  Let $N > 0$, and let $u_1$, $u_2$ be two solutions to the Schr\"odinger equation \eqref{schrodinger}, such that $u_j(t)$ has Fourier transform supported in the region $|\xi_j| \leq N$ for $j=1,2$.  Suppose also that the Fourier supports of $u_j(t)$ are separated by at least $\geq c N$.  Then for any $q > \frac{n+3}{n+1}$ we have the spacetime estimate
$$ \| u_1 u_2 \|_{L^q_{x,t}} \leq C(c) N^{n-\frac{n+2}{q}} P(u_1)^{1/2} P(u_2)^{1/2}.$$
\end{corollary}

\begin{proof} By scale invariance one can take $N=1$.  The claim then follows directly from Theorem \ref{main} and \eqref{prob}.
\end{proof}

Such a statement implies various bilinear estimates for $X^{s,b}$ norms for Schr\"odinger and wave equations, see e.g. \cite{tao:cone} for a discussion.  It is also likely that this sort of estimate has application to nonlinear Schr\"odinger equations; for instance, one can combine this estimate with the arguments in \cite{planchon} to obtain new well-posedness results for certain non-linear Schr\"odinger equations in Besov spaces.

Another application to Schr\"odinger equations was noted in \cite{tv:cone2}.  Indeed, from Theorem 2.1 in \cite{tv:cone2} and Theorem \ref{main} of this paper we see immediately that $H^s$ solutions to \eqref{schrodinger} converge pointwise to the initial data as $t \to 0$ for $n=2$ and $s > 2/5$; this improves upon the result of $s > 15/32$ given in that paper, but does not reach the conjectured level of $s \geq 1/4$.  In higher dimensions $n \geq 2$, a direct modification of the arguments in \cite{tv:cone2} gives convergence for $s > n/(n+3)$.

The author is a Clay Prize Fellow and is supported by the Packard Foundation.  The author also thanks Fabrice Planchon for helpful comments, and the anonymous referee for careful reading of the paper and many cogent suggestions (which have since been incorporated into the paper).

\section{Notation}\label{notation-sec}

If $X$ is a finite set, we use $\# X$ to denote its cardinality; if $X$ is a measurable set, we use $|X|$ to denote its Lebesgue measure.

If $(t,x) \in \R \times \R^n$ is a point in spacetime, we use $B((t,x),r)$ to denote the spacetime ball 
$$ B((t,x),r) := \{ (t',x') \in \R \times \R^n: |(t',x')-(t,x)| < r \}$$
and $D(x,r)$ to denote the spatial disk
$$ D(x,r) := \{ x' \in \R^n: |x'-x| < r \}.$$

We use $A \lesssim B$ or $A = O(B)$ to denote the estimate $|A| \leq CB$, where $C$ is a constant depending only on $n$.

Very shortly, our estimates shall involve a large parameter $R \gg 1$.
We shall use $A \lessapprox B$ to denote the estimate $A \leq C_\eps R^\eps B$
for all $\eps > 0$; in particular we note that $(\log R)^C \lessapprox 1$
for any $C$.

\section{Reduction to localized restriction estimates}\label{localized-sec}

We now begin the proof of Theorem \ref{main}.  Our arguments closely follow that of Wolff \cite{wolff:cone}, but with one additional twist near the end.  The argument is organized as follows.  In this section we make a preliminary reduction to the problem of obtaining sufficiently good localized restriction estimates, and then set up the induction argument we will use to obtain such estimates.  In the next section we recall the wave packet decomposition of Schr\"odinger waves, which has been fundamental to all of the recent developments in restriction theory for these waves.  Then, in Section \ref{inductive-sec}, we use the inductive hypothesis to strip away a certain ``localized'' component of the estimate, and reduce ourselves to considering only the ``global'' portion.  To estimate this global expression we perform a standard fine scale decomposition of space in Section \ref{fine-sec}, splitting the problem into obtaining a fine-scale estimate and then a coarse-scale estimate.  The fine-scale estimate is purely Fourier-analytic and is estimated using Plancherel's theorem in Section \ref{improv-sec}; our innovation here is to exploit an additional constraint on frequencies arising from the codimension 1 nature of the frequency space hypersurface $S$.  The coarse-scale estimate is a geometric
combinatorics estimate of Kakeya type, and is proven by the standard Bourgain-Wolff ``bush'' counting argument; the key point is that the constraint on frequencies from the fine-scale analysis translates\footnote{The linkage between fine scales and coarse scales is provided (heuristically, at least) by the \emph{dispersion relation}, which asserts that the frequency of a wave at fine scales determines the (group) velocity of that wave at coarse scales; in the physical interpretation of the Schr\"odinger equation, this relation is codified by de Broglie's law $p = \hbar \xi$.  To make this heuristic mathematically rigorous, the wave packet decomposition is an ideal tool.} to a constraint on directions in the coarse-scale estimate, thus restricting the bush to a hypersurface.  This puts us in the situation to apply Wolff's counting argument from \cite{wolff:cone}, which then concludes the proof.

We now turn to the details.
Fix $S_1$, $S_2$.  By a finite partition of $S_1$ and $S_2$, exploiting the compactness hypothesis, we may assume that $\diam(S_1), \diam(S_2) \ll \dist(S_1, S_2)$.  After a suitable rotation, scaling, and Gallilean transformation (the latter effectively translates $\xi$ by an arbitrary amount while keeping $S$ invariant), one may thus assume that
$$ S_1 := \{ (\tau, \xi) \in S: |\xi - e_1| \leq \frac{1}{100n} \}$$
and
$$ S_2 := \{ (\tau, \xi) \in S \leq \frac{1}{100n} \},$$
where $e_1$ is a standard unit basis vector.  We shall also need the slight enlargements
$$ \tilde S_1 := \{ (\tau, \xi) \in S: |\xi - e_1| \leq \frac{1}{50n} \}$$
and
$$ \tilde S_2 := \{ (\tau, \xi) \in S \leq \frac{1}{50n} \},$$

Following Wolff \cite{wolff:cone}, our first step is to reduce matters to proving a localized restriction estimate in which we are permitted to lose epsilon powers of the localization scale $R$.

\begin{definition}  We use $R^*_{S_1 \times S_2}(2 \times 2 \to q, \alpha)$ to denote the estimate
$$ \| \widehat{f_1 d\sigma_1} \widehat{f_2 d\sigma_2}\|_{L^q(B((t_0,x_0),R)} \leq C_{q,S_1,S_2,\alpha} R^\alpha \|f_1\|_{L^2(S_1, d\sigma_1)} \|f_2\|_{L^2(S_2, d\sigma_1)}$$
for all smooth $f_1$, $f_2$ on $S_1$, $S_2$, all $R \geq 1$, and all spacetime balls $B((t_0,x_0),R)$ of radius $R$.
\end{definition}

To prove Theorem \ref{main}, it suffices by standard ``epsilon-removal'' lemmas\footnote{For instance, one can apply Lemma 2.4 from \cite{tv:cone1}; see also Section 4 of \cite{borg:cone}, or Section 8 of \cite{krt}.  In all of these arguments (which are of Tomas-Stein type) the key fact is that the surface measures on $S_1$ and $S_2$ has a Fourier transform which decays at infinity; this is ultimately a consequence of the non-vanishing curvature of these surfaces.} to prove the local estimate
\be{alpha}
R^*_{S_1 \times S_2}(2 \times 2 \to \frac{n+3}{n+1}, \alpha).
\end{equation}
for all $\alpha > 0$.

To prove \eqref{alpha} we use Wolff's induction on scale argument.  It is easy to see that the above estimate must be true for sufficiently large $\alpha$; for instance, one can use the crude bound $\| \widehat{f_j d\sigma_j} \|_\infty \leq C \| f_j \|_2$ to obtain \eqref{alpha} for some large $\alpha$.  The claim will then follow (as in Wolff \cite{wolff:cone}) from the following inductive statement.

\begin{proposition}\label{inductive}  Suppose $\alpha > 0$ is such that \eqref{alpha} holds. Then we have
$$ R^*_{S_1 \times S_2}(2 \times 2 \to \frac{n+3}{n+1}, \max((1-\delta)\alpha, C\delta) + C \eps)$$
for all $0 < \delta, \eps \ll 1$, where the constants $C$ are independent of $\delta$ and $\eps$.
\end{proposition}

By choosing $\delta$ and $\eps$ suitably we may make $\max((1-\delta)\alpha, C\delta) + C \eps$ equal to $\alpha - c \alpha^2$ for some small absolute constant $c$.  Iterating this we thus see that the infimum of all $\alpha > 0$ for which \eqref{alpha} holds is zero, and the claim follows.

It remains to prove Proposition \ref{inductive}.  This will occupy the rest of the paper.

\section{The wave packet decomposition}\label{wave-sec}

As in the arguments\footnote{The basic idea of using wave packet decompositions to attack restriction and Bochner-Riesz type problems goes back to Fefferman and C\'ordoba.} of Bourgain \cite{borg:stein}, \cite{borg:cone}, Wolff \cite{wolff:cone}, and others, the next step is to decompose the functions $\widehat{f_1 d\sigma_1}$ and $\widehat{f_2 d\sigma_2}$ into wave packets concentrated on $R \times \sqrt{R}$ tubes.  

Fix $R \gg 1$ (the case $R \sim 1$ being trivial), and let $j=1,2$.   Let $\Z^n$ be the standard integer lattice in $\R^n$.  We shall need a spatial grid $X := R^{1/2} \Z^n$ and a velocity grid $V := R^{-1/2} \Z^n$.  We let $V_j$ be those velocities\footnote{Note that because of our normalization of the paraboloid \eqref{paraboloid}, the group velocity $v$ is exactly equal to the spatial frequency $\xi$; physically, this is just de Broglie's relation $mv = p = \hbar \xi$ under the normalization $m = \hbar = 1$.  Thus we will not bother to make much of a distinction between velocity and frequency in this argument.} $v \in V$ such that $(1, v)$ is normal to $\tilde S_j$, i.e. $(-\frac{1}{2}|v|^2, v) \in \tilde S_j$.

We shall work on the spacetime slab $[R/2,R] \times \R^n$.
We define a \emph{$\tilde S_j$-tube} to be any set of the form
$$ T := \{ (t,x): R/2 \leq t \leq R; |x - (x(T) + t v(T))| \leq R^{1/2} \},$$
where $x(T) \in X$ is the \emph{initial position} of $T$ and $v(T) \in V_j$ is the \emph{velocity}. 

We shall need the following standard wave packet decomposition (this is the parabola analogue of the cone decompositions in \cite{wolff:cone}, \cite{tao:cone}, \cite{krt}, and is also implicit in \cite{borg:stein}, \cite{tvv:bilinear}):

\begin{lemma}\label{decomp-lemma}  Let $j=1,2$, and let $f_j$ be a smooth function on $S_j$.  Then there exists a decomposition
\be{decomp}
\widehat{f_j d\sigma_j} = \sum_{T_j} c_{T_j} \phi_{T_j}
\end{equation}
where $T_j$ ranges over all $\tilde S_j$-tubes, the complex-valued co-efficients $c_{T_j}$ obey the $l^2$ bound
\be{l2}
(\sum_{T_j} |c_{T_j}|^2)^{1/2} \lesssim \|f_j\|_2,
\end{equation}
and for each $T_j$, the wave packets $\phi_{T_j}$ are free Schrodinger waves, where for each $R/2 \leq t \leq R$, the function $\phi_{T_j}(t)$ has Fourier transform supported on the set
\be{ft-local}
\{ \xi \in \R^n: \xi = v(T_j) + O(R^{-1/2}) \}
\end{equation}
(informally, $\phi_{T_j}$ has frequency $v(T_j) + O(R^{-1/2})$) and obeys the pointwise estimates
\be{decay-0}
|\phi_{T_j}(t,x)| \leq C_N R^{-n/4} (1 + \frac{|x - (x(T_j) + tv(T_j))|}{R^{1/2}})^{-N}
\end{equation}
for all $x \in \R^n$, and any $N > 0$.  In particular, outside of the tube
$$ R^\delta T_j := \{ (t,x): R/2 \leq t \leq R; |x - (x(T_j) + t v(T_j))| \leq R^{1/2+\delta} \},$$
we have the estimate
\be{decay-1}
|\phi_{T_j}(t,x)| \lessapprox R^{-100n}.
\end{equation}

Finally, any collection $\T_j$ of $\tilde S_j$-tubes, we have the probability estimate
\be{probability}
P(\sum_{T_j \in \T_j} \phi_{T_j}) \lesssim \# \T_j.
\end{equation}
\end{lemma}

\begin{proof}
We first prove the Lemma under the assumption that $f_j$ is supported on a cap of the form
\be{cap}
\{ (\tau,\xi) \in S_j: \xi = v + O(R^{-1/2}) \}
\end{equation}
for some fixed $v \in V_j$; this assumption will be removed at the end of this proof.

From the Poisson summation formula we may find a Schwartz function $\eta$ whose Fourier transform is supported in a disk $D(0,C) \subset \R^n$ such that $\sum_{k \in \Z^n} \eta(x-k) \equiv 1$.  Let $F(x) := \widehat{f_j d\sigma_j}(0,x)$ denote the initial data of $\widehat{f_j d\sigma_j}$.
We thus have the decomposition
$$ F(x) = \sum_{x_0 \in X_j} \eta(\frac{x-x_0}{R^{1/2}}) F(x).$$
Observe that the spatial Fourier transform of $\eta(\frac{x-x_0}{R^{1/2}}) F(x)$ is supported on a disk $\{ \xi \in \R^n: \xi = v + O(R^{-1/2})\}$.  Thus if we let $u_{x_0}$ be the unique Schr\"odinger wave with initial data $u_{x_0}(0,x) := \eta(\frac{x-x_0}{R^{1/2}}) \widehat{f_j d\sigma_j}(0,x)$, then we have the decomposition
$$ \widehat{f_j d\sigma_j}(t,x) = \sum_{x_0 \in X_j} u_{x_0}(t,x).$$
Now let $T_j$ be a $\tilde S_j$-tube with $v(T_j) = v$.  We write $c_{T_j} := R^{n/4} M F(x(T_j))$, where 
$$ MF(x) := \sup_{r>0} \frac{1}{|D(x,r)|} \int_{D(x,r)} |F|$$
is the Hardy-Littlewood maximal  function of $F$, and write $\phi_{T_j} := u_{x(T_j)}/c_{T_j}$.  Thus we have
$$ \widehat{f_j d\sigma_j} = \sum_{T_j: v(T_j) = v} c_{T_j} \phi_{T_j},$$
thus giving a decomposition \eqref{decomp} (setting $c_{T_j} = \phi_{T_j} = 0$ for $v(T_j) \neq v$).  Since $F$ has Fourier transform supported in the disk \eqref{ft-local}, it enjoys a reproducing formula of the form $F = F * \psi$ where the reproducing kernel $\psi = \psi_v$ has Fourier support in a (slight enlargement of) the disk \eqref{ft-local}, and obeys the pointwise bounds 
$$|\psi(x)| \leq C_N R^{-n/2} (1 + |x|/R^{1/2})^{-N}$$
for any $N \geq 0$.  From this it is easy to see that $MF(x) \sim MF(x')$ whenever $|x-x'| \lesssim R^{1/2}$.  Thus
$$ \sum_{T_j: v(T_j) = v} |c_{T_j}|^2 \lesssim \int |MF(x)|^2\ dx \lesssim \|F\|_2^2
\lesssim \|f_j\|_2^2$$
by the Hardy-Littlewood maximal  inequality and Plancherel's theorem; this gives \eqref{l2}.

By construction, the Fourier transform of $\phi_{T_j}(0)$ (and hence $\phi_{T_j}(t)$ for any $t$) is supported in the set \eqref{ft-local}.  Now we prove \eqref{decay-0}.  By construction, it suffices to show the pointwise estimate
\be{dc}
|u_{x_0}(t,x)| \leq C_N (1 + \frac{|x - (x_0 + tv)|}{R^{1/2}})^{-N} M F(x_0)
\end{equation}
for all $x_0 \in X$ and $t \sim R$.  By translation invariance we may take $x_0 = 0$.

There are several ways to prove this estimate; for instance, one can observe that \eqref{cap} is contained in an $O(R^{-1}) \times O(R^{-1/2})$ disk with normal $(1,v)$ and use some form of the uncertainty principle.  Another way to argue is as follows.
From the fundamental solution of the free Schr\"odinger equation we have an integral representation of the form
$$ u_{0}(t,x) = C t^{-n/2}
\int e^{i C |x-y|^2/t} \eta(\frac{y}{R^{1/2}}) F(y)\ dy.$$
Recall the reproducing formula $u_0 = u_0 * \psi_v$.  Thus we have
$$ u_{0}(t,x) = C t^{-n/2}
\int K_v(x-y) \eta(\frac{y}{R^{1/2}}) F(y)\ dy$$
where $K_v$ is the kernel
$$ K_v(x) := \int \int e^{i C |x-y|^2/t} \phi_v(y)\ dy.$$
A routine stationary phase computation\footnote{The reader may wish to simplify the calculation by first taking advantage of Gallilean invariance to reduce to the case $v=0$, and then using the scale invariance of the Schr\"odinger equation to reduce to the case $R=1$.}, using the decay and Fourier support properties of $\phi_v$, gives the bounds
$$ |K_v(x)| \leq C_N R^{-n/2} (1 + \frac{x-vt}{R^{1/2}})^{-N}$$
for all $N \geq 0$.  The claim then follows \eqref{dc} from a direct computation.

The estimate \eqref{decay-1} follows from \eqref{decay-0}, so it remains to prove \eqref{probability}.  Since the probability is time-invariant, it suffices to show that
$$ \int |\sum_{T_j \in \T_j: v(T_j) = v} \phi_{T_j}(0,x)|^2\ dx \lesssim \# \T_j.$$
But this follows directly from \eqref{decay-0}, since the tubes $T_j$ with fixed velocity $v(T_j) = v$ all have distinct initial positions $x(T_j)$, which are separated by $\gtrsim R^{1/2}$.

Now we remove the hypothesis that $f_j$ was supported in a cap \eqref{cap}.  For general $f_j$, we may of course decompose $f_j = \sum_{v \in V_j} f_{j,v}$, where each $f_{j,v}$ is supported in the cap \eqref{cap} associated to $v$, and we have the $L^2$ bound
\be{fff}
\sum_{v \in V_j} \| f_{j,v}\|_2^2 \sim \| f_j\|_2^2.
\end{equation}
One can then apply the previous arguments to $f_{j,v}$, obtaining a decomposition
$$ f_{j,v} = \sum_{T_j: v(T_j) = v} c_{T_j} \phi_{T_j}$$
obeying all the above properties.  Summing over all $v$ we obtain a decomposition \eqref{decomp} of $f_j$, which then  obeys \eqref{l2} thanks to \eqref{fff}.  The properties \eqref{decay-0}, \eqref{decay-1}, and the Fourier support in \eqref{ft-local} have all been proven, so it remains to show \eqref{probability}.  But we have already proven the special case
$$
P(\sum_{T_j \in \T_j: v(T_j) = v} \phi_{T_j}) \lesssim \# \{ T_j \in \T_j: v(T_j) = v \}$$
for all $v \in V_j$; the claim then follows by summing in $v$ and exploiting the frequency space orthogonality (via the support property \eqref{ft-local}).
\end{proof}

We can now begin the proof of Proposition \ref{inductive} in earnest.  Fix $\alpha$, and let $Q_R$ denote the cylinder
$$ Q_R := \{ (t,x): R/2 \leq t \leq R; |x| \leq R \}.$$ 
It will suffice to prove the estimate
$$
\| \widehat{f_1 d\sigma_1} \widehat{f_2 d\sigma_2} \|_{L^{\frac{n+3}{n+1}}(Q_R)} \lessapprox (R^{(1-\delta)\alpha} + R^{C \delta}) 
\| f_1 \|_2 \|f_2\|_2
$$
for all smooth $f_1$, $f_2$ on $S_1$, $S_2$, since any ball of radius $R$ can be covered by $O(1)$ translates of $Q_R$.  Here and in the sequel our
implicit constants in $\lessapprox$ or $\lesssim$ are allowed to depend on  
$\delta$.

Fix $f_1$, $f_2$; we may normalize $\|f_1\|_2 = \|f_2\|_2 = 1$.  We apply Lemma \ref{decomp-lemma} to both $f_1$ and $f_2$, writing
$$ f_j = \sum_{T_j} c_{T_j} \phi_{T_j}$$
for $j=1,2$, where $T_j$ ranges over $\tilde S_j$-tubes.  It thus suffices to show that
$$
\|\sum_{T_1} \sum_{T_2} c_{T_1} c_{T_2} \phi_{T_1} \phi_{T_2} \|_{L^{\frac{n+3}{n+1}}(Q_R)}
\lessapprox R^{(1-\delta)\alpha} + R^{C \delta}.
$$
We first remove some minor portions of this sum.  Let us first consider the contribution when $T_1$ and $T_2$ are both disjoint from $B(0, CR)$.  In this case the bound \eqref{decay-0} gives bounds of $O(R^{-100n})$ for both $\phi_{T_1}$ and $\phi_{T_2}$, with the bound improving even more as $T_1$ and $T_2$ move away from $B(0,CR)$.  Since the coefficients $c_{T_1}$, $c_{T_2}$ are bounded by \eqref{l2}, the contribution of this case is easily seen to be acceptable.

A similar argument disposes of the case where $T_1$ is disjoint from $B(0,CR)$ and $T_2$ intersects $B(0,CR)$, as in this case there are only $O(R^{2n})$ possible values of $T_2$.  Similarly when $T_2$ is disjoint from $B(0,CR)$ and $T_1$ intersects $B(0,CR)$.  Thus we may henceforth restrict ourselves to tubes which intersect $B(0,CR)$.  In particular, the number of tubes $T_1$ under consideration is now only $O(R^{2n})$, and similarly for $T_2$.

We can now eliminate the contribution of the terms where $c_{T_1} = O(R^{-100n})$ or $c_{T_2} = O(R^{-100n})$, since those terms can be easily controlled just by using $L^\infty$ bounds on $\phi_{T_1}$, $\phi_{T_2}$ (from e.g. \eqref{decay-0}).  Thus we only need to restrict ourselves to the tubes $T_1$ where $R^{-100n} \lesssim c_{T_1}\lesssim 1$, and similarly for $T_2$.

By pigeonholing the interval $[R^{-100n}, 1]$ dyadically into $O(\log R)$ groups, and noting that $\log R \approx 1$, we may thus restrict the $T_1$ summation to the tubes where $c_{T_1} \sim \gamma_1$ for some fixed $R^{-100n} \lesssim \gamma_1 \lesssim 1$.  Let $\T_1$ denote the set of all tubes $T_1$ of this form; from \eqref{l2} we have $(\# \T_1)^{1/2} \lesssim \gamma_1^{-1}$.  We may as well assume that $c_{T_1} = \gamma_1$ for these tubes $T_1 \in \T_1$, since we can absorb the factor $c_{T_1}/\gamma_1$ harmlessly into $\phi_{T_1}$.  Similarly, we may restrict the tubes $T_2$ to a collection $\T_2$ with $(\# \T_2)^{1/2} \lesssim \gamma_2^{-1}$ and $c_{T_2} = \gamma_2$ for all $T_2 \in \T_2$, for some
$R^{-100n} \lesssim \gamma_2 \lesssim 1$.  It thus suffices to prove 

\begin{proposition}\label{p-prop}  We have the estimate
\be{p-sum}
\| \sum_{T_1 \in \T_1} \sum_{T_2 \in \T_2} \phi_{T_1} \phi_{T_2} \|_{L^{\frac{n+3}{n+1}}(Q_R)} \lessapprox (R^{(1-\delta)\alpha} + R^{C \delta}) 
(\# \T_1)^{1/2} (\# \T_2)^{1/2}
\end{equation}
for all collections $\T_1$, $\T_2$ of $\tilde S_1$-tubes and $\tilde S_2$-tubes respectively, such that all the tubes intersect $B(0,CR)$.
\end{proposition}

It remains to prove this Proposition.  This will be done in the next few sections.

\section{Localization of tubes, and the inductive argument}\label{inductive-sec}

We now utilize the inductive hypothesis \eqref{alpha}.  The idea (due to Wolff \cite{wolff:cone}) is to give each wave packet $\phi_{T_1}$ and $\phi_{T_2}$ a slightly smaller ball of radius $R^{1-\delta}$ which it can ``exclude'' via the inductive hypothesis; it will then suffice to verify the $L^p$ estimate on the exterior of these balls.  This is similar to the ``two-ends'' reduction used in the Kakeya problem, see e.g. \cite{wolff:kakeya}.

We turn to the details.  We may cover the cylinder $Q_R$ by about $O(R^{C\delta})$ finitely overlapping spacetime balls $B$ of radius $R^{1-\delta}$; let $\B$ denote the collection of such balls.  We can thus estimate the left-hand side of \eqref{p-sum} extremely crudely\footnote{Clearly we may improve on this by replacing the $l^1$ summation over balls $B$ with an $l^{(n+3)/(n+1)}$ summation.  This refinement is exploited in the endpoint theory, see \cite{tao:cone}, but is unnecessary for the non-endpoint case.} by
\be{p-sum-2}
\sum_{B \in \B} \| \sum_{T_1 \in \T_1} \sum_{T_2 \in \T_2} \phi_{T_1} \phi_{T_2} \|_{L^{\frac{n+3}{n+1}}(B)}. 
\end{equation}

Suppose we have some relation $\sim$ between the tubes in $\T_1 \cup \T_2$ and balls in $\B$; we will specify this relation much later in the argument, but roughly we will associate $T \sim B$ if the contribution of $\phi_T$ to the bilinear expression $\sum_{T_1 \in \T_1} \sum_{T_2 \in \T_2} \phi_{T_1} \phi_{T_2}$ is ``concentrated'' in $B$.  We can then estimate \eqref{p-sum-2} by the ``local part''

\be{p-sum-local}
\sum_{B \in \B} \| (\sum_{T_1 \in \T_1: T_1 \sim B} \phi_{T_1}) (\sum_{T_2 \in \T_2: T_2 \sim B} \phi_{T_2}) \|_{L^{\frac{n+3}{n+1}}(B)}
\end{equation}
and the ``global part''
\be{p-sum-global}
\sum_{B \in \B} \| \sum_{T_1 \in \T_1, T_2 \in \T_2: T_1 \not \sim B \hbox{ or } T_2 \not \sim B} \phi_{T_1} \phi_{T_2} \|_{L^{\frac{n+3}{n+1}}(B)}.
\end{equation}
Consider the contribution of the local portion \eqref{p-sum-local}.  From the probability estimate \eqref{probability} we see that for each $B \in \B$ and $j=1,2$, $\sum_{T_j \in \T_j: T_j \sim B} \phi_{T_j}$ is a free Schr\"odinger wave with probability
$$ P(\sum_{T_j \in \T_j: T_j \sim B} \phi_{T_j}) \lesssim \#\{T_j \in \T_j: T_j \sim B \}.$$
By applying the induction hypothesis \eqref{alpha}, we may thus bound 
\eqref{p-sum-local} by
$$ \eqref{p-sum-local} \lessapprox
\sum_{B \in \B} R^{(1-\delta)\alpha} (\#\{T_1 \in \T_1: T_1 \sim B \})^{1/2}
(\#\{T_2 \in \T_2: T_2 \sim B \})^{1/2},$$
which by Cauchy-Schwarz becomes
$$ \eqref{p-sum-local} \lessapprox R^{(1-\delta)\alpha}
(\sum_{B \in \B} \sum_{T_1 \in \T_1: T_1 \sim B} 1)^{1/2}
(\sum_{B \in \B} \sum_{T_2 \in \T_2: T_2 \sim B} 1)^{1/2}.$$
Thus, if we make 

\begin{assumption}\label{ass} For all $T \in \T_1 \cup \T_2$, we have
\be{sim-bound}
\# \{ B \in \B: T \sim B \} \lessapprox 1,
\end{equation}
\end{assumption}

then we can bound \eqref{p-sum-local} by
$$ \eqref{p-sum-local} \lessapprox R^{(1-\delta)\alpha} (\# \T_1)^{1/2} (\# \T_2)^{1/2}$$
which is acceptable.

Roughly speaking, Assumption \ref{ass} asserts that each tube $T \in \T_1 \cup \T_2$ is allowed to exclude $\lessapprox 1$ balls $B$ from the summation in \eqref{p-sum-global}.  It is thus natural to select $\sim$ so that each tube $T$ excludes the ball $B$ in which its ``contribution'' to \eqref{p-sum} is ``greatest''; this will become clearer when we define $\sim$ in Section \ref{kakeya-sec}.

It remains to estimate \eqref{p-sum-global}.  It will suffice to show that
\be{p-piece}
\| \sum_{T_1 \in \T_1, T_2 \in \T_2: T_1 \not \sim B \hbox{ or } T_2 \not \sim B} \phi_{T_1} \phi_{T_2} \|_{L^{\frac{n+3}{n+1}}(B)} \lessapprox
R^{C\delta} (\# \T_1)^{1/2} (\# \T_2)^{1/2}
\end{equation}
for all $B \in \B$, since the claim then follows by summing in $B$.  Note that $\alpha$ no longer plays any role; we will not need the induction hypothesis \eqref{alpha} in the remainder of the argument.  Also, we can now freely lose powers of $R^\delta$ in what follows.

Fix $B$; it remains to prove \eqref{p-piece}. By the triangle inequality, it will suffice to prove that
\be{p-piece-a}
\| \sum_{T_1 \in \T_1: T_1 \not \sim B} \sum_{T_2 \in \T_2} \phi_{T_1} \phi_{T_2} \|_{L^{\frac{n+3}{n+1}}(B)} \lessapprox R^{C\delta}
(\# \T_1)^{1/2} (\# \T_2)^{1/2}
\end{equation}
and
\be{p-piece-b}
\| \sum_{T_2 \in \T_2: T_2 \not \sim B} \sum_{T_1 \in \T_1: T_1 \sim B} \phi_{T_1} \phi_{T_2} \|_{L^{\frac{n+3}{n+1}}(B)} \lessapprox
R^{C\delta} (\# \T_1)^{1/2} (\# \T_2)^{1/2}.
\end{equation}
The two claims are proven similarly (the expression \eqref{p-piece-b} is slightly smaller, but the extra constraint $T_1 \sim B$ turns out to play no significant role), and so we will content ourselves with proving \eqref{p-piece-a}.
(The definition of the equivalence relation $\sim$ will be symmetric with respect to $\T_1$ and $\T_2$).

We follow Wolff's strategy of obtaining the bilinear $L^{\frac{n+3}{n+1}}$ estimate by interpolating between bilinear $L^1$ and $L^2$ estimates.  The bilinear $L^1$ estimate follows easily from linear $L^2$ estimates:

\begin{lemma}\label{l1-easy}  We have
$$
\| \sum_{T_1 \in \T_1: T_1 \not \sim B} \sum_{T_2 \in \T_2} \phi_{T_1} \phi_{T_2} \|_{L^1(B)} \lessapprox R
(\# \T_1)^{1/2} (\# \T_2)^{1/2}.$$
\end{lemma}

\begin{proof}
By H\"older's inequality it suffices to show that
$$ 
\| \sum_{T_1 \in \T_1: T_1 \not \sim B} \phi_{T_1} \|_{L^2(B)} \lessapprox R^{1/2} (\# \T_1)^{1/2}$$
and
$$ 
\| \sum_{T_2 \in \T_2} \phi_{T_2} \|_{L^2(B)} \lessapprox R^{1/2} (\# \T_2)^{1/2}.$$
But these follow directly from \eqref{probability} and an integration in time (since $B$ is contained in the slab $[-R,R] \times \R^n$). 
\end{proof}

From Lemma \ref{l1-easy} and H\"older's inequality (or the log-convexity of $L^p$ norms), it will suffice to prove the $L^2$ estimate 
\be{prod}
\| \sum_{T_1 \in \T_1: T_1 \not \sim B} \sum_{T_2 \in \T_2} \phi_{T_1} \phi_{T_2} \|_{L^2(B)} \lessapprox R^{C\delta} R^{-(n-1)/4}
(\# \T_1)^{1/2} (\# \T_2)^{1/2};
\end{equation}
note how this uses the choice of exponent $\frac{n+3}{n+1}$.

The exponent $R^{-(n-1)/4}$ is best possible.  To see this, let $\pi$ denote the spacetime disk
$$ \pi := \{ (t, x_1 e_1): t, x_1 = O(R) \},$$
and consider the example when $\T_1$ consists of the $O(\sqrt{R})$ tubes with velocity $e_1$ which intersect $\pi$, while $\T_2$ similarly consists of the $O(\sqrt{R})$ tubes with velocity $-e_1$ which also intersect the $\pi$ plane. By \eqref{decay-0}, the left hand side is essentially of magnitude $O(R^{-n/2})$ on a $O(\sqrt{R})$-neighbourhood $\pi$ (which thus has volume $R^{(n+3)/2}$), and the numerology of \eqref{prod} follows.  (This is of course the same counterexample which shows that the exponent $\frac{n+3}{n+1}$ is best possible; see \cite{tvv:bilinear}, \cite{wolff:cone}).

\section{Fine-scale decomposition}\label{fine-sec}

In the previous part of the argument, we have decomposed the cylinder $Q_R$ (which is essentially a spacetime ball of radius $R$) into slightly smaller balls $B$ of radius $R^{1-\delta}$ in order to utilize the induction hypothesis.  To continue the argument we must decompose $B$ into much smaller balls, namely balls of radius $\sqrt{R}$, to fully exploit the spatial localization of the tubes $T$.  Specifically, we cover (a slight dilate of) $Q_R$ by a finitely overlapping collection $\q$ of balls of radius $\sqrt{R}$.  Squaring \eqref{prod}, it thus suffices to show that
\be{2-bound}
\sum_{q \in \q: q \subset 2B} \| \sum_{T_1 \in \T_1: T_1 \not \sim B} \sum_{T_2 \in \T_2} \phi_{T_1} \phi_{T_2} \|_{L^2(q)}^2 \lessapprox R^{C\delta} R^{-(n-1)/2}
(\# \T_1) (\# \T_2).
\end{equation}

First consider the contribution to \eqref{2-bound} of the case where $T_1 \cap R^\delta q = \emptyset$.  In this case, it is easy to see from \eqref{decay-1} and the triangle inequality that this contribution is certainly acceptable. Thus we only need to consider the terms in \eqref{2-bound} where $T_1$ intersects $R^\delta q$.  Similarly we only need to consider the terms where $T_2$ intersects $R^\delta q$.  

It remains to show
\be{2-bound-local}
\sum_{q \in \q: q \subset 2B } \| \sum_{T_1 \in \T^{\not \sim B}_1(q)} \sum_{T_2 \in \T_2(q)} \phi_{T_1} \phi_{T_2} \|_{L^2(q)}^2 \lessapprox R^{C\delta} R^{-(n-1)/2}
(\# \T_1) (\# \T_2)
\end{equation}
where
\bas
 \T_j(q) &:= \{ T_j \in \T_j : T_j \cap R^\delta q \neq \emptyset \} \hbox{ for } j=1,2\\
 \T^{\not \sim B}_1(q) &:= \{ T_1 \in \T_1(q) : T_1 \not \sim B \}.
\end{align*}

We now do some dyadic pigeonholing, first on the multiplicity of the tubes $T_1$, $T_2$ through $q$, and then on the multiplicity of the balls $q$ within $T_1$.  For any dyadic numbers\footnote{By \emph{dyadic number} we mean an integer power of two.} $1 \leq \mu_1, \mu_2 \lessapprox R^{100n}$, let $\q(\mu_1, \mu_2) \subset \q$ denote the set
$$ \q(\mu_1, \mu_2) := \{ q \in \q: \mu_j \leq \# \T_j(q) < 2\mu_j \hbox{ for } j=1,2 \},$$
thus the $\q(\mu_1, \mu_2)$ cover all the balls $q \in \q$ for which the summand in \eqref{2-bound-local} is non-zero.  Since there are only $O(\log R)^2 \approx 1$ possible values of $(\mu_1, \mu_2)$, it thus suffices to show that
$$
\sum_{q \in \q(\mu_1, \mu_2): q \subseteq 2B} \| \sum_{T_1 \in \T^{\not \sim B}_1(q)} \sum_{T_2 \in \T_2(q)} \phi_{T_1} \phi_{T_2} \|_{L^2(q)}^2 \lessapprox R^{C\delta} R^{-(n-1)/2}
(\# \T_1) (\# \T_2)
$$
for all $\mu_1, \mu_2$.

Fix $\mu_1, \mu_2$.  For any $T_1 \in \T_1$, let $\lambda(T_1,\mu_1,\mu_2) $ denote the integer
$$ \lambda(T_1, \mu_1, \mu_2) := \# \{ q \in \q(\mu_1, \mu_2): T_1 \cap R^\delta q \neq \emptyset \},$$
and for every dyadic number $1 \leq \lambda_1 \leq R^{100n}$, let $\T_1[\lambda_1,\mu_1,\mu_2]$ denote the set
\be{tlmm-def}
\T_1[\lambda_1,\mu_1,\mu_2] := \{ T_1 \in \T_1: \lambda_1 \leq \lambda(T_1,\mu_1,\mu_2) < 2 \lambda_1 \}.
\end{equation}
Since there are only $O(\log R) \approx 1$ values of $\lambda_1$, it thus suffices to show that
\be{pigeonholed}
\sum_{q \in \q(\mu_1, \mu_2): q \subset 2B} \| \sum_{T_1 \in \T^{\not \sim B}_1(q) \cap \T_1[\lambda_1,\mu_1,\mu_2]} \sum_{T_2 \in \T_2(q)} \phi_{T_1} \phi_{T_2} \|_{L^2(q)}^2 \lessapprox R^{C\delta} R^{-(n-1)/2}
(\# \T_1) (\# \T_2)
\end{equation}
for all $\lambda_1$.  (We could also pigeonhole the multiplicity of balls in $T_2$ in a similar manner, but this will turn out to be unnecessary).

Fix $\lambda_1$.  We still have to prove \eqref{pigeonholed}.  At this point we pause to recall how the analogous argument of Wolff \cite{wolff:cone} proceeded for the cone (for which the tubes $T_1$, $T_2$ are constrained to point in null directions).  Firstly, by a Plancherel argument (similar to an argument of Mockenhaupt \cite{mock:cone}), Wolff observed the local estimate
\be{wolff-local}
\begin{split}
\| \sum_{T_1 \in \T^{\not \sim B}_1(q) \cap \T_1[\lambda_1,\mu_1,\mu_2]} &\sum_{T_2 \in \T_2(q)} \phi_{T_1} \phi_{T_2} \|_{L^2(q)}^2 \lessapprox \\
&R^{C\delta} R^{-(n-1)/2}
\#(T^{\not \sim B}_1(q) \cap T_1[\lambda_1,\mu_1,\mu_2])^2 (\# \T_2(q));
\end{split}
\end{equation}
this is basically a consequence of the fact that for fixed $T_1$, the functions $\phi_{T_1} \phi_{T_2}$ are almost orthogonal on $q$.  From \eqref{wolff-local}
it would then suffice to show the combinatorial estimate
\be{chunk}
\sum_{q \in \q(\mu_1, \mu_2): q \subset 2B} \#(\T^{\not \sim B}_1(q) \cap \T_1[\lambda_1,\mu_1,\mu_2])^2 (\# \T_2(q))\lessapprox R^{C\delta}
(\# \T_1) (\# \T_2).
\end{equation}
This estimate is true in the case of the cone (see the Remarks section) but does not appear to hold for the paraboloid case.  To resolve this difficulty we need to sharpen the local estimate \eqref{wolff-local}; this is the purpose of the next section.

\section{An improved local estimate}\label{improv-sec}

Before we present our improved version of the local estimate \eqref{wolff-local}, let us begin with an informal discussion.  Suppose we wish to estimate a quantity of the form
\be{22}
\| \sum_{T_1 \in \T_1} \sum_{T_2 \in \T_2} \phi_{T_1} \phi_{T_2} \|_{L^2_{t,x}}^2,
\end{equation}
where we shall be careless about exactly what region of spacetime we are integrating over.  We can expand this expression as
$$ \sum_{T_1 \in \T_1} \sum_{T_2 \in \T_2} \sum_{T'_1 \in \T_1} \sum_{T'_2} \langle \phi_{T_1} \phi_{T_2}, \phi_{T'_1} \phi_{T'_2} \rangle.$$
Now if $T_1$ has velocity $\xi_1$, then the spacetime Fourier transform $\phi_{T_1}$ should be supported near the point $(\xi_1, -\frac{1}{2}|\xi_1|^2)$ in $S$.  Similarly if $T_2$ has velocity $\xi_2$, $T'_1$ has velocity $\xi'_1$, and $T_2$ has velocity $\xi'_2$.  From Parseval's formula, we thus expect the above inner product to be very small unless $\xi_1 + \xi_2$ is close to $\xi'_1 + \xi'_2$ and $|\xi_1|^2 + |\xi_2|^2$ is close to $|\xi'_1|^2 + |\xi'_2|^2$.

Suppose we fix two of the frequencies, say $\xi_1$ and $\xi'_1$.  Then the relation $\xi_1 + \xi_2 = \xi'_1 + \xi'_2$ will correlate $\xi_2$ and $\xi'_2$, in the sense that either of these two frequencies will determine the other.  This basic observation is already enough to give a bound for \eqref{22} which is proportional to $(\# \T_1)^2 (\# \T_2)$, and by making these ideas slightly more rigorous one can soon obtain the bound \eqref{wolff-local}.  However, as we will soon see, we can do better by also exploiting the additional constraint $|\xi_1|^2 + |\xi_2|^2 = |\xi'_1|^2 + |\xi'_2|^2$ to remove one more degree of freedom on the collection $\T_2$, which will eventually make this collection behave sufficiently similar to the collection of tubes in a light cone that Wolff's argument will apply.

We need some notation.  Let $\Omega_1, \Omega_2 \subseteq \R^n$ denote the spatial frequency regions
\bas
\Omega_1 &:= \{ \xi \in \R^n: |\xi - e_1| \leq \frac{1}{20n} \}\\
\Omega_2 &:= \{ \xi \in \R^n: |\xi + e_1| \leq \frac{1}{20n} \};
\end{align*}
note these are slightly larger than the spatial frequency supports of $\tilde S_1$ and $\tilde S_2$ respectively.  For any $\xi_1 \in \Omega_1, \xi_2 \in \Omega_2$, let $\pi(\xi_1, \xi'_2) \subseteq \Omega_1$ denote the set
\be{pi-def}
\begin{split}
\pi(\xi_1, \xi'_2) := \{ &\xi'_1 \in \Omega_1: \xi_1 + \xi_2 = \xi'_1 + \xi'_2, |\xi_1|^2 + |\xi_2|^2 = |\xi'_1|^2 + |\xi'_2|^2 \\
&\hbox{ for some } \xi_2 \in \Omega_2 \};
\end{split}
\end{equation}
one can interpret this set as being equivalent to the set of all parallelograms with two vertices in (a slight enlargement of) $\tilde S_1$ and two vertices in (a slight enlargement of) $\tilde S_2$.

A little algebra shows that $\pi(\xi_1, \xi'_2)$ is contained in the $n-1$-dimensional hyperplane in $\R^n$ which contains $\xi_1$ and is orthogonal to $ \xi'_2 - \xi_1$ (cf. the calculations in \cite{borg:schrodinger}, \cite{vargas:restrict}, \cite{vargas:2}, \cite{tvv:bilinear}), or in other words \footnote{This orthogonality is not absolutely essential to the argument; what is important (particularly in the proof of Lemma \ref{t2}) is that the set $\pi(\xi_1, \xi'_2)$ is contained in a hypersurface which is transverse to $\xi'_2 - \xi_1$, or indeed to any vector in $\Omega_2 - \Omega_1$.}
\be{pi-orthog}
\langle \xi'_1 - \xi_1, \xi'_2 - \xi_1 \rangle_{\R^n} = 0 \hbox{ whenever } \xi'_1 \in \pi(\xi_1, \xi'_2).
\end{equation}
Indeed, the points $\xi_1, \xi'_1, \xi_2, \xi'_2$ form a rectangle in $\R^n$.

For any ball $q \in \q(\mu_1, \mu_2)$ and any two frequencies $\xi_1 \in \Omega_1$, $\xi'_2 \in \Omega_2$, let $\T^{\not \sim B}_1(q, \lambda_1, \mu_1, \mu_2, \xi_1, \xi'_2)$ denote the collection of those tubes $T_1 \in \T^{\not \sim B}_1(q) \cap \T_1[\lambda_1, \mu_1, \mu_2]$ such that the velocity $v(T_1)$ of $T_1$ is within $O(R^{C\delta} R^{-1/2})$ of the set $\pi(\xi_1, \xi_2')$.  Let $\nu(q,\lambda_1, \mu_1, \mu_2)$ denote the quantity
\be{nu-def}
\nu(q, \lambda_1, \mu_1, \mu_2) := \sup_{\xi_1 \in \Omega_1; \xi'_2 \in \Omega_2} \# \T^{\not \sim B}_1(q, \lambda_1, \mu_1, \mu_2, \xi_1, \xi'_2).
\end{equation}

We now prove the following refinement of \eqref{wolff-local}.

\begin{lemma}\label{terry-local}  For any $q \in \q(\mu_1, \mu_2)$, we have
\bas
&\left\| \sum_{T_1 \in \T^{\not \sim B}_1(q) \cap \T_1[\lambda_1,\mu_1,\mu_2]} \sum_{T_2 \in \T_2(q)} \phi_{T_1} \phi_{T_2} \right\|_{L^2(q)}^2 
\lessapprox \\
&R^{C\delta} R^{-(n-1)/2}
\nu(q,\lambda_1, \mu_1, \mu_2) (\#(\T^{\not \sim B}_1(q) \cap \T_1[\lambda_1,\mu_1,\mu_2])) (\# \T_2(q)).
\end{align*}
\end{lemma}

\begin{proof}  Our arguments here have certain similarities to those in \cite{mock:cone}, \cite{borg:schrodinger}, \cite{vargas:restrict}, \cite{vargas:2}, \cite{tvv:bilinear}, and can ultimately be traced back to the $L^4$ theory of Fefferman, Sj\"olin, and C\'ordoba.

For brevity, let us write 
\bas
\T'_1 &:= \T^{\not \sim B}_1(q) \cap \T_1[\lambda_1,\mu_1,\mu_2]\\
\T'_2 &:= T_2(q)\\
\nu &:= \nu(q,\lambda_1, \mu_1, \mu_2).
\end{align*}
Since the $L^2(q)$ norm is bounded by the global $L^2_{t,x}$ norm, it suffices to show that
\be{phi-4}
\| \sum_{T_1 \in \T'_1} \sum_{T_2 \in \T'_2} \phi_{T_1} \phi_{T_2} \|_{L^2_{t,x}}^2 \lessapprox
R^{C\delta} R^{-(n-1)/2} (\# \T'_1) (\# \T'_2) \nu.
\end{equation}
Note that a straightforward calculation using Plancherel's theorem shows that $\phi_{T_1} \phi_{T_2}$ is globally in $L^2_{t,x}$. Indeed, from \eqref{decay-0} we see that the spacetime Fourier transform of $\phi_{T_1}$ is of the form $\widehat{f_{T_1}\ d\sigma_1}$, where $f_{T_1}$ is supported on a cap $\{ (\tau,\xi) \in S: \xi = v(T_1) + O(R^{-1/2}) \}$ and has magnitude $O(R^{n/4})$.  Similarly for $\phi_{T_2}$.  A computation using the transversality of $S_1$ and $S_2$ thus shows that the spacetime Fourier transform of $\phi_{T_1} \phi_{T_2}$ is supported on the spacetime ball
\be{fs}
B( (-\frac{1}{2} |v(T_1)|^2, v(T_1)) + (-\frac{1}{2} |v(T_2)|^2, v(T_2)), C R^{-1/2} )
\end{equation}
and has magnitude $O(R^{1/2})$.  In particular we have
\be{phiphi}
\| \phi_{T_1} \phi_{T_2} \|_{L^2_{t,x}} \lesssim R^{-(n-1)/4}.
\end{equation}
We now return to \eqref{phi-4}.  We expand the left-hand side as
\be{phi-2-phi}
|\sum_{T_1, T'_1 \in \T'_1} \sum_{T_2, T'_2 \in \T'_2} 
\langle \phi_{T_1} \phi_{T_2}, \phi_{T'_1} \phi_{T'_2} \rangle_{L^2_{t,x} }|.
\end{equation}
From \eqref{phiphi} and Cauchy-Schwarz we see that the inner product is $O(R^{-(n-1)/2})$.  On the other hand, from the Fourier support \eqref{fs}, we see that the inner product vanishes unless 
\be{vvv}
v(T_1) + v(T_2) = v(T'_1) + v(T'_2) + O(R^{-1/2})
\end{equation}
and
$$ |v(T_1)|^2 + |v(T_2)|^2 = |v(T'_1)|^2 + |v(T'_2)|^2 + O(R^{-1/2}).$$
In particular, we see (using the separation of $\Omega_1$ and $\Omega_2$)
that for fixed $T_1$, $T'_2$, the velocity $v(T'_1)$ must lie within $O(R^{-1/2})$ of the hyperplane $\pi(v(T_1), v(T'_2))$.  In particular for fixed $T_1$, $T'_2$ there are at most $O(\nu)$ choices for $v(T'_1)$, and hence $O(R^{C\delta} \nu)$ choices of $T'_1$ (since by construction of $\T'_1$, $T'_1$ must intersect $R^\delta q$).    For fixed $T_1$, $T'_2$, $T_2$ there are at most $O(1)$ choices of $T'_1$ by \eqref{vvv}, and hence $O(R^{C\delta})$ choices of $T'_1$.  Combining all these facts together, we see that we can bound \eqref{phi-2-phi} by
$$ (\# \T'_1) (\# \T'_2) R^{C\delta} \nu R^{C\delta} R^{-(n-1)/2},$$
as desired.
\end{proof}

To conclude the proof of Theorem \ref{main}, it thus remains to prove the combinatorial (Kakeya-type) estimate
\be{combinatorial}
\begin{split}
\sum_{q \in \q(\mu_1,\mu_2): q \subset 2B} &\nu(q,\lambda_1, \mu_1, \mu_2) (\#(\T^{\not \sim B}_1(q) \cap \T_1[\lambda_1,\mu_1,\mu_2])) (\# \T_2(q))\\
&\lessapprox R^{C\delta} (\# \T_1) (\# \T_2)
\end{split}
\end{equation}
for an appropriate choice of equivalence relation $\sim$ obeying Assumption
\ref{ass}.
This will occupy the next section.

\section{The combinatorial estimate}\label{kakeya-sec}

We now prove the combinatorial estimate \eqref{combinatorial}.  Interestingly, this estimate is of a comparable level of difficulty to the corresponding combinatorial estimate\footnote{In our notation, the combinatorial estimate in \cite{wolff:cone} is essentially \eqref{chunk}, but with the tubes $T_j$ restricted to light rays.  See also the remarks section.} in \cite{wolff:cone}, and in particular does not need any additional Kakeya-type information.  (The numerology is similar to the $(n+2)/2$ Kakeya estimate in \cite{wolff:kakeya}, but the argument here seems simpler than the ``hairbrush'' argument in \cite{wolff:kakeya}, though of a somewhat similar flavor).

We first need to define the relation $\sim$.  For each tube $T_1 \in \T_1[\lambda_1, \mu_1, \mu_2]$, let $B(T_1,\lambda_1,\mu_1,\mu_2)$ be the ball in $\B$ which maximizes the quantity 
$$ \# \{ q \in \q(\mu_1, \mu_2): T_1 \cap R^\delta q \neq \emptyset; q \cap B(T_1, \lambda_1, \mu_1, \mu_2) \neq \emptyset \}.$$
From the pigeonhole principle and \eqref{tlmm-def}, we observe that
\be{q-large}
 \# \{ q \in \q(\mu_1, \mu_2): T_1 \cap R^\delta q \neq \emptyset; q \cap B(T_1, \lambda_1, \mu_1, \mu_2) \neq \emptyset \} \gtrapprox R^{-C\delta} \lambda_1.
\end{equation}
We define the relation $\sim_{\lambda_1, \mu_1, \mu_2}$ between tubes in $\T_1$ and balls in $\B$ by defining $T_1 \sim_{\lambda_1,\mu_1,\mu_2} B'$ if $T_1 \in \T_1[\lambda_1, \mu_1, \mu_2]$ and $B' \subseteq 10 B(T_1,\lambda_1,\mu_1,\mu_2)$; note that this definition is independent of the ball $B$ which appeared in the previous section.  Clearly for each tube $T_1$ there are at most $O(1)$ balls $B'$ such that $T_1 \sim_{\lambda_1, \mu_1, \mu_2} B'$.  Then we define $T_1 \sim B'$ if one has $T_1 \sim_{\lambda_1, \mu_1, \mu_2} B'$ for some dyadic $\lambda_1, \mu_1, \mu_2$; it is then clear that \eqref{sim-bound} holds for $T \in \T_1$.  We then define $\sim$ between $\T_2$ and $\B$ by a completely symmetrical procedure (although we will not need $\sim$ for $\T_2$ here as we are proving \eqref{p-piece-a} instead of \eqref{p-piece-b}).  

Now we prove \eqref{combinatorial}.
By definition of $\q(\mu_1,\mu_2)$, we have 
\be{t2-card}
\# \T_2(q) \lessapprox \mu_2
\end{equation}
 for all $q$ in \eqref{combinatorial}.  Also, by Fubini's theorem and \eqref{tlmm-def}, we have
\be{fubini}
\begin{split}
\sum_{q \in \q(\mu_1, \mu_2): q \subset 2B} \#(\T^{\not \sim B}_1(q) \cap \T_1[\lambda_1])
&\leq \sum_{q \in \q(\mu_1, \mu_2)} \#(\T_1(q) \cap \T_1[\lambda_1])\\
&= \sum_{T_1 \in \T_1[\lambda_1]} \# \{ q \in \q(\mu_1, \mu_2): T_1 \cap R^\delta q \neq \emptyset \} \\
&\lessapprox \sum_{T_1 \in \T_1[\lambda_1,\mu_1,\mu_2]} \lambda_1 \\
&\leq (\# \T_1) \lambda_1.
\end{split}
\end{equation}
Thus to prove \eqref{combinatorial} it will suffice to show that
\be{nu-mult}
\nu(q_0,\lambda_1,\mu_1,\mu_2) \lessapprox R^{C\delta} \frac{\# \T_2}{\lambda_1 \mu_2}
\end{equation}
for all $q_0 \in \q(\mu_1, \mu_2)$ with $q_0 \subset 2B$.

It remains to prove \eqref{nu-mult}, which we shall do using a ``bush'' argument centered at $q_0$.  Fix $q_0 \in \q(\mu_1, \mu_2)$ with $q_0 \subset 2B$, and let $\xi_1 \in \Omega_1$, $\xi'_2 \in \Omega_2$ be arbitrary.  
Let $\T'_1$ denote the set
\be{tp1-def}
\T'_1 := \T^{\not \sim B}_1(q_0, \lambda_1, \mu_1, \mu_2, \xi_1, \xi'_2)
\end{equation}
defined in Section \ref{improv-sec}.  By \eqref{nu-def}, it suffices to show that
\be{t-b}
\# \T'_1 \lessapprox R^{C\delta} \frac{\# \T_2}{\lambda_1 \mu_2}.
\end{equation}

Let $T_1 \in \T'_1$.  By construction, we have $T_1 \in \T^{\not\sim B}_1(q_0)$ and $T_1 \in \T_1[\lambda_1, \mu_1, \mu_2]$.  In particular, we have $T_1 \cap R^\delta q_0 \neq \emptyset$, and $B \not \subset 10 B(T_1, \lambda_1, \mu_1, \mu_2)$.  In particular, since $q_0 \subset 2B$, we have
$$ \dist(q_0, 2B(T_1, \lambda_1, \mu_1, \mu_2)) \gtrapprox R^{-C\delta} R.$$
By \eqref{q-large}, we thus have
$$ \# \{ q \in \q(\mu_1, \mu_2): T_1 \cap R^\delta q \neq \emptyset; \dist(q_0, q) \gtrapprox R^{-C\delta} R \} \gtrapprox R^{-C\delta} \lambda_1.$$
On the other hand, by the definition of $\q(\mu_1, \mu_2)$, for each $q \in \q(\mu_1, \mu_2)$ there are $\gtrapprox \mu_2$ tubes $T_2$ in $\T_2$ which intersect $R^\delta q$.  Thus we have
$$ \# \{ (q, T_2) \in \q(\mu_1, \mu_2) \times \T_2: T_1 \cap R^\delta q, T_2 \cap R^\delta q \neq \emptyset; \dist(q_0, q) \gtrapprox R^{-C\delta} R \} \gtrapprox R^{-C\delta} \lambda_1 \mu_2.$$
Summing over all $T_1$ in $\T'_1$, we obtain
\be{com}
\# \{ (q, T_1, T_2) \in \q \times \T'_1 \times \T_2: T_1 \cap R^\delta q, T_2 \cap R^\delta q \neq \emptyset; \dist(q_0, q) \gtrapprox R^{-C\delta} R \} \gtrapprox R^{-C\delta} \lambda_1 \mu_2 \# \T'_1.
\end{equation}
Now we make the following crucial geometric observation, which is analogous to the geometric observation used in \cite{wolff:cone} that a light ray can transversally intersect a light cone in at most one point:

\begin{lemma}\label{t2}  For each $T_2 \in \T_2$, we have
$$ \# \{ (q, T_1) \in \q \times \T'_1: T_1 \cap R^\delta q, T_2 \cap R^\delta q \neq \emptyset; \dist(q_0, q) \gtrapprox R^{-C\delta} R \} \lessapprox R^{C\delta}.$$
\end{lemma}

\begin{proof}
Let $(t_0,x_0)$ and $(t,x)$ denote the centers of $q_0$ and $q$ respectively.  Since $T_1$ intersects both $R^\delta q_0$ and $R^\delta q$, and $\dist(q_0, q) \gtrapprox R^{-C\delta} R$, we see that
$$ R^{-C\delta} R \lessapprox |t-t_0| \lessapprox R$$
and
$$ x - x_0 = v(T_1)(t-t_0) + O(R^{C\delta} R^{1/2}).$$
On the other hand, since $T_1 \in \T'_1$, we see from \eqref{tp1-def} that $v(T_1)$ lies within $O(R^{C\delta} R^{-1/2})$ of $\pi(\xi_1, \xi'_2)$.  Thus we have
$$ \dist( \frac{x-x_0}{t-t_0}, \pi(\xi_1, \xi'_2)) \lessapprox R^{C\delta} R^{-1/2}.$$
On the other hand, if we let $e := \xi'_2 - \xi_1$, then from \eqref{pi-orthog} we see that
$$\langle \xi'_1 - \xi_1, e \rangle_{\R^n} = 0 \hbox{ for all } \xi'_1 \in \pi(\xi_1, \xi_2),$$
and hence
$$ \langle \frac{x-x_0}{t-t_0} - \xi_1, e \rangle_{\R^n} \lessapprox R^{C\delta} R^{-1/2}.$$
We may rearrange this as
$$ \langle (t-t_0,x-x_0), (-\langle \xi_1,e \rangle_{\R^n}, e) \rangle_{\R \times \R^n} \lessapprox R^{C\delta} R^{1/2}.$$
Thus $(t,x)$ lies within $O(R^{C\delta} R^{1/2})$ of the $n$-dimensional hyperplane $\Pi$ in $\R \times \R^{n-1}$ which passes through $(t_0,x_0)$ and which is normal to $(-\langle \xi_1,e\rangle_{\R^n}, e)$.  But since $\xi_1 \in \Omega_1$, $\xi'_2 \in \Omega_2$, we see that $e$ is within $1/5n$ of $-2e_1$, and $-\langle \xi_1, e \rangle$ is within $1/5n$ of $+2$.  Since $v(T_2)$ is within $1/5n$ of $-e_1$, we thus see that $T_2$ makes an angle of $\sim 1$ with respect to $\Pi$.  Since $\dist((t,x), T_2) \lessapprox R^{C\delta} R^{1/2}$, we thus see that $(t,x)$ is thus constrained to lie within a ball of radius $R^{C\delta} R^{1/2}$.  This means that there are only at most $O(R^{C\delta})$ choices for $q$.  For each fixed $q$ there are at most $O(R^{C\delta})$ choices for $T_1$, and the claim follows. 
\end{proof}

Combining this Lemma with \eqref{com} we see that
$$ R^{C\delta} \# \T_2 \gtrapprox R^{-C\delta} \lambda_1 \mu_2 \# \T'_1$$
and \eqref{t-b} follows.  This concludes the proof of Theorem \ref{main}.
\endprf

\section{Remarks}\label{remarks-sec}

\begin{itemize}

\item The proof of Theorem \ref{main} is very similar to the argument in \cite{wolff:cone}.  Indeed, one can compare the arguments as follows.  For the cone, the passage to localized restriction estimates, wave packet decomposition, induction on scales, and fine scale decomposition works almost exactly the same as with the parabola, the only major difference being that the tubes are now
oriented along light rays\footnote{Also, the tubes have a more interesting internal structure, being composed of somewhat thinner $1 \times R^{1/2} \times R$ ``plates'', but this ends up not being very relevant to the argument which follows.  See \cite{wolff:cone}, \cite{tao:cone} for further discussion.}.  For the localized estimate, \eqref{wolff-local} is used instead of Lemma \ref{terry-local}.  This requires us to prove \eqref{chunk}.  Using \eqref{fubini} and \eqref{t2-card} as in Section \ref{kakeya-sec}, one reduces to showing that
$$\#(\T^{\not \sim B}_1(q_0) \cap \T_1[\lambda_1,\mu_1,\mu_2]) \lessapprox R^{C\delta} \frac{\# \T_2}{\lambda_1 \mu_2}.$$
Arguing as in Section \ref{kakeya-sec}, this reduces to showing the estimate
$$ \# \{ (q, T_1) \in \q \times \T_1: T_1 \cap R^\delta q, T_2 \cap R^\delta q \neq \emptyset; \dist(q_0, q) \gtrapprox R^{-C\delta} R \} \lessapprox R^{C\delta}.$$

(compare with Lemma \ref{t2}).  But this follows in the cone case since the tubes $T_1$ which intersect $R^\delta q$ are contained in a $R^{1/2+\delta}$-neighborhood of a light cone; since the tube $T_2$ is concentrated around a light ray, intersects $T_1$ transversally and at a distance $\gtrapprox R^{-C\delta} R$ from the vertex of this light cone, the claim then follows from elementary geometry. 

\item It may well be possible to eliminate much of the pigeonholing in the above argument, and perhaps even eradicate the epsilon loss in Theorem \ref{main}.  (See for instance \cite{tao:informal} for a non-pigeonholed version of the argument for the cone in \cite{wolff:cone}, and \cite{tao:cone} for the endpoint result).  However, it seems difficult to access the $\nu$ parameter without this pigeonholing, and we do not know how to remove the epsilons in the paraboloid case.

\item The geometric properties of the paraboloid which were used in the above argument (and especially in Lemma \ref{t2}) are easily seen to be robust under small perturbations of the paraboloid.  In particular, one can easily obtain Theorem \ref{main} for all disjoint compact subsets of a compact hypersurface of \emph{elliptic type} as defined in \cite{vargas:2}, \cite{tvv:bilinear}, providing that the parameter $\eps$ used to define elliptic type is sufficiently small. We sketch this as follows.  Let $S$ be a surface of elliptic type; after some linear transformations, this means that $S$ is of the form
$$ S := \{ (\tau,\xi) \in \R \times \R^n: \tau = -\frac{1}{2} |\xi|^2 + \eps f(\xi)\}$$
where the error function $f(\xi)$ is smooth, and $\eps$ is a sufficiently small parameter (depending on the smooth norms of $f$ and on the size and separation of $S_1$, $S_2$).  In other words, $S$ is a small perturbation of the paraboloid \eqref{paraboloid}.  This means that the dispersion relation $v = h(\xi)$ between the group velocity $v$ and the frequency $\xi$ is not quite the identity (in fact, it is given by $h(\xi) := \xi - \eps \nabla f(\xi)$), but it will still be a homeomorphism and a small perturbation of the identity on $S_1 \cup S_2$ if $\eps$ is small enough.
Aside from making this distinction between velocity and frequency, the arguments in Sections \ref{localized-sec}-\ref{fine-sec} are essentially unchanged.  In Section \ref{improv-sec}, the set $\pi(\xi_1, \xi'_2)$ must be replaced by
\bas \pi_S(\xi_1,\xi'_2) :=
& \{ \xi'_1 \in \Omega_1: \xi_1 + \xi_2 = \xi'_1 + \xi'_2, |\xi_1|^2 + |\xi_2|^2 = |\xi'_1|^2 + |\xi'_2|^2 \\
+ &2 \eps(f(\xi_1)+f(\xi_2)-f(\xi'_1)-f(\xi'_2))
\hbox{ for some } \xi_2 \in \Omega_2 \},
\end{align*}
but this is easily seen to be a small smooth perturbation of $\pi(\xi_1, \xi'_2)$.  Actually, because the dispersion relation $v = h(\xi)$ is no longer the identity, the relevant set is not $\pi_S(\xi_1,\xi'_2)$ but rather $h(\pi_S(\xi_1,\xi'_2))$, but this is still a small smooth perturbation of $\pi(\xi_1,\xi'_2)$, and in particular retains the key property of lying in a hypersurface transverse to $\Omega_2 - \Omega_1$.

Now the remainder of the argument continues as before, with the obvious modifications, until we reach Lemma \ref{t2}.  Now $(t,x)$ will not lie within $O(R^{C\delta} R^{1/2})$ of a hyperplane in spacetime, but instead it will lie
within $O(R^{C\delta} R^{1/2})$ of a conic manifold\footnote{In the special case when $S$ is a sphere, then this conic manifold is in fact a circular cone, although the aperture and orientation of this cone depends on $\xi_1$ and $\xi'_2$.} consisting of the union of the lines through vertex $(t_0, x_0)$ which have velocity in $h(\pi_S(\xi_1,\xi'_2))$.  If $\eps$ is sufficiently small, this manifold is still transverse to $T_2$, and the remainder of the argument proceeds as before.

\item Once we have the above bilinear restriction theorems for arbitrary disjoint compact subsets of surfaces of elliptic type, we can use the machinery of \cite{tvv:bilinear} to derive the analogue of Corollary \ref{rest-conj} for all compact hypersurfaces of elliptic type.  After some finite partitions of unity and some affine linear transformations, we may thus obtain Corollary \ref{rest-conj} for all compact surfaces for which all the principal curvatures strictly positive.  In particular, the restriction conjecture for the sphere $S^n$ in $\R \times \R^n$ is true for all $q > 2(n+3)/(n+1)$.  

\item It is also extremely likely that the same argument works when some of the principal curvatures are strictly negative; indeed, by combining this argument with the argument for the cone, it seems plausible that one should be able to obtain good restriction estimates for all surfaces in which at most one principal curvature vanishes at any given point.  In particular, one should be able to obtain bilinear restriction theorems for all non-degenerate conic sections when $q > 2(n+3)/(n+1)$ (thus providing a bilinear analogue of the linear theory in \cite{strichartz:restrictionquadratic}).  If so, this would likely give near-optimal bilinear $L^p$ null form estimates for the wave equation (see \cite{tao:cone} for a discussion).

\item It seems likely that these arguments also give some new progress on the Bochner-Riesz problem for paraboloids and spheres (see e.g. \cite{borg:stein} for a discussion), but we have not pursued this question. 

\end{itemize}

\end{document}